\documentclass[11pt]{amsart}
\usepackage{amssymb, amsmath, graphicx}

\newtheorem{thm}{Theorem}[section]

\newtheorem{lem}[thm]{Lemma}
\newtheorem{prop}[thm]{Proposition}

\theoremstyle{definition}
\newtheorem{defi}[thm]{Definition}

\theoremstyle{remark}

\numberwithin{equation}{section}

\newcommand{\al}{\alpha}

\newcommand{\De}{\Delta}
\newcommand{\Ga}{\Gamma}
\newcommand{\Om}{\Omega}
\newcommand{\om}{\omega}
\renewcommand{\th}{\theta}
\newcommand{\ze}{\zeta}

\newcommand{\s}{\mathbf s}
\newcommand{\n}{\mathbf n}
\newcommand{\bZ}{\mathbf Z}

\newcommand{\LL}{\mathcal L}
\newcommand{\C}{\mathbb C}

\newcommand{\Z}{\mathbb Z}

\newcommand{\R}{\mathbb R}
\newcommand{\CP}{{\mathbb C}{\mathbb P}}
\newcommand{\del}{\partial}

\newcommand{\x}{\times}
\newcommand{\co}{\colon\thinspace}

\DeclareMathOperator{\rk}{rk}
\newcommand{\Spin}{\mathrm Spin}

\begin{document}

\title{On Lens Spaces and Their Symplectic Fillings}

\author{Paolo Lisca}
\address{Dipartimento di Matematica\\ Universit\`a di Pisa \\I-56127
Pisa, ITALY} 
\email{lisca@dm.unipi.it}

\thanks{The author is a member of EDGE, Research Training Network
HPRN-CT-2000-00101, supported by The European Human Potential
Programme. The author's research was partially supported by MURST}

\subjclass{Primary 57R17; Secondary 53D35}
\date{February 25, 2002}

\begin{abstract} 
The standard contact structure $\xi_0$ on the three--sphere $S^3$ is
invariant under the action of $\Z/p\Z$ yielding the lens space
$L(p,q)$, therefore every lens space carries a natural quotient
contact structure $\overline\xi_0$. A theorem of Eliashberg and McDuff
classifies the symplectic fillings of ($L(p,1),{\overline\xi_0})$ up
to diffeomorphism. Here we announce a generalization of that result
to every lens space. In particular, we give an explicit handlebody
decomposition of every symplectic filling of
$(L(p,q),{\overline\xi_0})$ for every $p$ and $q$. Our results imply:
\begin{enumerate}
\item[(a)]
There exist infinitely many lens spaces $L(p,q)$ with $q\not= 1$ such
that $(L(p,q),{\overline\xi_0})$ admits only one symplectic filling up
to blowup and diffeomorphism.
\item[(b)]
For any natural number $N$, there exist infinitely many lens spaces
$L(p,q)$ such that $(L(p,q),{\overline\xi_0})$ admits more than $N$
symplectic fillings up to blowup and diffeomorphism.
\end{enumerate}
\end{abstract}

\maketitle

\section{Introduction}
\label{s:intro}

Four--dimensional symplectic fillings are objects of central interest
in symplectic and contact topology. They can be used to prove
tightness of contact structures on three--dimensional
manifolds~\cite{El3}, to distinguish tight contact
structures~\cite{LM}, and they arise in symplectic cut--and--paste
constructions~\cite{Sy}. Symplectic fillings also have an intrinsic
interest. For example, they behave like closed symplectic
four--manifolds from the point of view of the Seiberg--Witten
invariants and this fact has a number of consequences~\cite{KM, Li1,
Li2}.

A difficult and interesting problem in this area is the diffeomorphism
classification of the symplectic fillings of a given contact
three--manifold. Eliashberg and McDuff solved this problem in some
cases~\cite{El2, McD1}. The purpose of this note is to announce a
generalization of their results.

A (coorientable) \emph{contact three--manifold} is pair $(Y,\xi)$,
where $Y$ is a three--manifold and $\xi\subset TY$ is a
two--dimensional distribution given as the kernel of a one--form
$\al\in\Om^1(Y)$ such that $\al\wedge d\al$ is a volume form. 

A \emph{symplectic filling} of a closed contact three--manifold $(Y,\xi)$
is pair $(X,\om)$ consisting of a smooth, connected four--manifold $X$
with $\del X=Y$ and a symplectic form $\om$ on $X$ such that $\om|_\xi
\not= 0$ at every point of $\del X$. Moreover, if $X$ is oriented by
$\om\wedge\om$ and $\xi=\{\al=0\}$, the boundary orientation on $Y$
must coincide with the orientation induced by $\al\wedge d\al$ (which
only depends on $\xi$). 

A basic fact to bear in mind is that any diffeomorphism classification
of symplectic fillings will always be ``up to blowups''. This is
because a blowup, i.e.~a connected sum with $\overline\CP^2$, is a
local operation which can be performed in the symplectic
category. Therefore, if $(X,\om)$ is a symplectic filling of
$(Y,\xi)$ then $\widehat X = X\# N\CP^2$, $N\geq 1$, carries a
symplectic form $\widehat\om$ such that $(\widehat X, \widehat\om)$ is
still a symplectic filling of $(Y,\xi)$.

The \emph{standard contact structure} $\xi_0$ on the three--sphere is
the $2$--dimensional distribution of complex lines tangent to
$S^3\subset\C^2$. The unit four--ball $B^4\subset\C^2$ endowed with
the restriction of the standard K\"ahler form on $\C^2$ is a
symplectic filling of $(S^3,\xi_0)$. The following result, due to
Eliashberg, yields the diffeomorphism classification of the symplectic
fillings of $(S^3,\xi_0)$. 
\begin{thm}[\cite{El2}]
Let $(X,\om)$ be a symplectic filling of $(S^3,\xi_0)$. Then, $X$ 
is diffeomorphic to a blowup of $B^4$.
\end{thm}
The standard contact structure $\xi_0$ is invariant under the natural
action of $U(2)$ on $S^3$. Thus, $\xi_0$ is a fortiori invariant under
the induced action of the subgroup
\begin{equation*}
G_{p,q}=
\{\left(\begin{smallmatrix}
\ze & 0\\
0 & \ze^q
\end{smallmatrix}\right)
\ |\ \ze^p=1\}\subset U(2),
\end{equation*}
where $p, q\in\Z$. It follows that when $p> q\geq 1$ and $p, q$ are
coprime, $\xi_0$ induces a contact structure ${\overline\xi_0}$ on the
lens space $L(p,q)=S^3/G_{p,q}$. 

Let $D_p$ denote the disk bundle over the $2$--sphere with Euler class
$p$. It is not difficult to construct a symplectic form $\om$ on
$D_{-p}$ such that $(D_{-p},\om)$ is a symplectic filling of
$(L(p,1),{\overline\xi_0})$. Another symplectic filling of
$(L(4,1),{\overline\xi_0})$ is given by ${\mathcal
C}=\CP^2\setminus\nu(C)$, where $\nu(C)$ is a pseudo--concave
neighborhood of a smooth conic $C\subset\CP^2$, endowed with the
restriction of the standard K\"ahler form. The four--manifold
${\mathcal C}$ is easily shown to have the handlebody decomposition
given by the right--hand side of Figure~\ref{f:example1} in
section~\ref{s:examples}. 

The following theorem is a by--product of Dusa McDuff's classification
of closed symplectic four--manifolds containing symplectic spheres of
non--negative self--intersection. It yields the diffeomorphism
classification of the symplectic fillings of
$(L(p,q),{\overline\xi_0})$.
\begin{thm}[\cite{McD1}]\label{t:mcduff-quoted}
Let $(X,\om)$ be a symplectic filling of $(L(p,1),{\overline\xi_0})$.
Then, $X$ is diffeomorphic to a blowup of:
\begin{enumerate}
\item[(a)]
$D_{-p}$ if $p\not= 4$,
\item[(b)]
$D_{-4}$ or ${\mathcal C}$ if $p=4$.
\end{enumerate}
\end{thm}

\section{Symplectic fillings of $(L(p,q),{\overline\xi_0})$}
\label{s:results}

From now on, we shall assume that $p$ and $q$ are two coprime integers
such that
\begin{equation}\label{e:assumptions}
p>q\geq 1,\qquad
\frac p{p-q} = b_1 - \cfrac{1}{b_2 -
       \cfrac{1}{\ddots -
        \cfrac{1}{b_k}
}},
\end{equation}
where $b_i\geq 2$, $i=1,\ldots, k$. Note that there is a unique
continued fraction expansion of $\frac p{p-q}$ as
in~\eqref{e:assumptions}. We shall use the shorthand $[b_1,\ldots,
b_k]$ for such continued fractions. 

\begin{figure}[ht]
\includegraphics[height=1.2in]{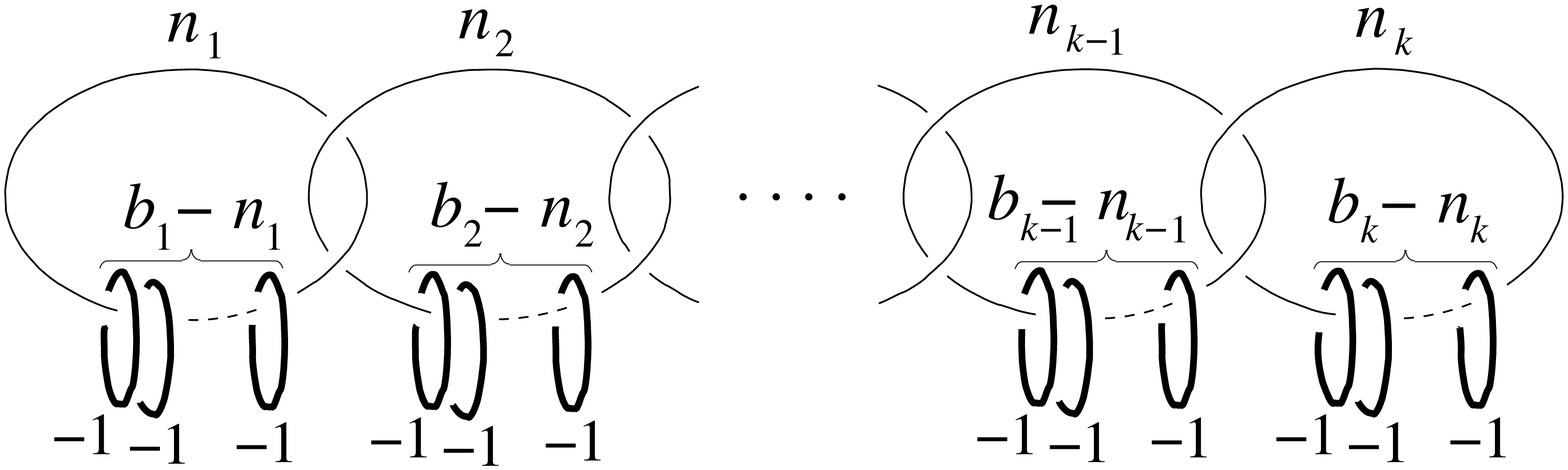}
\caption{The manifold $W_{p,q}(\n)$}
\label{f:W} 
\end{figure}

\medskip
Let $\n=(n_1,\ldots, n_k)$ be a $k$--tuple of non--negative integers
such that $[n_1,\ldots, n_k]=0$. We can view the ``thin'' framed link
in Figure~\ref{f:W} as a three--dimensional surgery presentation of a
closed, oriented three--manifold $N(\n)$. The assumption $[n_1,\ldots,
n_k]=0$ ensures the existence of an orientation--preserving
diffeomorphism
\begin{equation}\label{e:varphi}
\varphi\co N(n_1,\ldots, n_k)\to S^1\x S^2. 
\end{equation}

\medskip
\begin{defi}\label{d:thefillings}\label{d:W}
Let $W_{p,q}(\n)$ be the smooth four--manifold with boundary obtained
by attaching $2$--handles to $S^1\x D^3$ along the framed link
$\varphi({\mathbf L})\subset S^1\x S^2$, where ${\mathbf L}\subset
N(\n)$ is the thick framed link in Figure~\ref{f:W}.
\end{defi}

Observe that the manifold $W_{p,q}(\n)$ does not depend
on the choice of the diffeomorphism $\varphi$, because every
self--diffeomorphism of $S^1\x S^2$ extends to $S^1\x D^3$~\cite{Gl}.
Also, since ${\mathbf L}$ is canonically embedded in $S^3$, the
framings of its components can be identified with integers.

\medskip
Define $\bZ_{p,q}\subset\Z^k$ as follows:
\[
\bZ_{p,q}=\{(n_1,\ldots,n_k)\in\Z^k\ |\ [n_1,\ldots,n_k]=0, 
\  0\leq n_i\leq b_i,\ i=1,\ldots, k\}.
\]

We are now ready to state our main result.

\begin{thm}[\cite{Li3}]\label{t:main}
\begin{enumerate}
\item[(a)] Let $\n\in \bZ_{p,q}$. Then, the four--manifold $W_{p,q}(\n)$
carries a symplectic form $\om$ such that $(W_{p,q}(\n),\om)$ is a
symplectic filling of $(L(p,q), {\overline\xi_0})$. Moreover,
$W_{p,q}(\n)$ does not contain smoothly embedded spheres with
self--intersection $-1$.
\item[(b)] If $(W,\om)$ is a symplectic filling of $(L(p,q),
{\overline\xi_0})$, then $W$ is diffeomorphic to a blowup of
$W_{p,q}(\n)$ for some $\n\in \bZ_{p,q}$.
\end{enumerate}
\end{thm}

It is interesting to contrast Theorem~\ref{t:main} with certain
properties of the versal deformations of certain isolated complex
bidimensional singularities. The lens space $L(p,q)$, endowed with the
contact structure ${\overline\xi_0}$, can be viewed as the link of a
cyclic quotient singularity of type $(p,q)$. Stevens~\cite{St} used
deep results of Koll\'ar and Shepherd-Barron~\cite{KSB} to show that
the irreducible components of the reduced base space $S_{\text{red}}$
of the versal deformation of such a singularity are in one--to--one
correspondence with the set $\bZ_{p,q}$ (see also~\cite{Ch}).

Every irreducible component of $S_{\text{red}}$ gives a smoothing of
the singularity and therefore a symplectic filling of
$(L(p,q),{\overline\xi}_0)$. It follows by Theorem~\ref{t:main} that
every smoothing of a quotient $(p,q)$--singularity is diffeomorphic to
a manifold of the form $W_{p,q}(\n)$.

Indeed, it seems very likely that the smoothing coming from the
irreducible component indexed by $\n\in \bZ_{p,q}$ be diffeomorphic to
$W_{p,q}(\n)$. For example, it is easy to check that $\bZ_{p,q}$
always contains the $k$--tuple $(1,2,\ldots,2,1)$ as long as $k >
1$. The corresponding irreducible component of $S_{\text{red}}$ is
called the \emph{Artin component}. The associated smoothing is
diffeomorphic to the canonical resolution the singularity $R_{p,q}$,
which is a well--known four--manifold obtained by plumbing according
to a certain linear graph (see for example~\cite{BPV}). It is not
difficult to see that $R_{p,q}$ is orientation--preserving
diffeomorphic to $W_{p,q}(1,2,\ldots,2,1)$.

\medskip
If $\frac pq=[a_1,\ldots,a_h]$ with $a_i\geq 5$, $i=1,\ldots,k$, it
turns out that $\bZ_{p,q}=\{(1,2,\ldots,2,1)\}$. Thus, in this case
$S_{\text{red}}$ coincides with the Artin component, a fact which had
been conjectured by Koll\'ar~\cite{Ko}.

Under the same assumption, Theorem~\ref{t:main} implies that every
symplectic filling of $(L(p,q),{\overline\xi_0})$ is diffeomorphic to
a blowup of $W_{p,q}(1,2,\ldots,2,1)$. In particular, we see that
\emph{there exist infinitely many lens spaces $L(p,q)$ with $q\not= 1$ such
that $(L(p,q),{\overline\xi_0})$ admits only one symplectic filling up
to blowup and diffeomorphism}.

\medskip
On the other hand, suppose that $k\geq 4$, $b_2,\ldots,b_{k-2}\geq 3$
and $b_k\geq k-2$. Then,
\[
\n(r,s)=(1,\overbrace{2,\ldots,2}^r,3,\overbrace{2,\ldots,2}^s,1,s+2)
\in \bZ_{p,q}
\]
for every $0\leq r,s\leq k-4$ with $r+s+4=k$. Moreover, one can 
easily check that
\begin{equation}\label{e:rank}
\rk H_2\left(W_{p,q}(\n(r,s));\Z\right)
 = \sum_{i=1}^k b_i - 2k -s.
\end{equation}
Since the manifolds $W_{p,q}(\n)$ do not contain exceptional spheres
and therefore are not blowups of one another, Equation~\eqref{e:rank}
implies that $(L(p,q),{\overline\xi_0})$ admits at least $k-3$
symplectic fillings up to blowup and diffeomorphism. This shows that
\emph{for any natural number $N$, there exist infinitely many lens
spaces $L(p,q)$ such that $(L(p,q),{\overline\xi_0})$ admits more than
$N$ symplectic fillings up to blowup and diffeomorphism}.

\section{Examples}\label{s:examples}

In this section we illustrate Theorem~\ref{t:main} by analyzing its
implications in a few particular cases.

\medskip
Given a $k$--tuple of positive integers
\[
(n_1,\ldots, n_{s-1}, n_s, n_{s-1},\ldots, n_k)
\]
with $n_s=1$, we say that the $k-1$--tuple 
\[
(n_1,\ldots, n_{s-1}-1, n_{s+1}-1,\ldots, n_k)
\]
is obtained by a \emph{blowdown at $n_s$} (with the obvious meaning of
the notation when $s=1$ or $s=k$). The reverse process is a
\emph{blowup}.

We say that a $k$--tuple of positive integers $(n_1,\ldots, n_k)$ is
\emph{admissible} if the continued fraction $[n_1,\ldots, n_k]$ makes
sense, i.e.~if none of the denominators appearing in $[n_1,\ldots,
n_k]$ vanishes.

\medskip
The following Lemma~\ref{l:0-sequence} can be proved by an easy
induction.

\begin{lem}\label{l:0-sequence}
Let $(n_1,\ldots, n_k)$ be an admissible $k$--tuple of positive
integers. Then, $[n_1,\ldots, n_k]=0$ if and only if $(n_1,\ldots,
n_k)$ is obtained from $(0)$ by a sequence of blowups.
\end{lem}

\noindent\textbf{First Example.}  Let us apply
Theorem~\ref{t:main} to determine all the symplectic fillings of
$(L(p,1),{\overline\xi_0})$. We have
\[
\frac p{p-1}=[\overbrace{2,\ldots,2}^{p-1}].
\]
Using Lemma~\ref{l:0-sequence} one can easily check that:
\[
\bZ_{p,1}=
\begin{cases}
\{(1,2,\ldots,2,1)\},\ \text{if $p\not=4$},\\
\{(1,2,1),(2,1,2)\},\ \text{if $p=4$}.
\end{cases}
\]
Therefore, Theorem~\ref{t:main} implies that the symplectic fillings
of $(L(p,1),{\overline\xi_0})$ up to blowup and diffeomorphism are
$W_{p,1}(1,2,\ldots,2,1)$ for $p\not=4$ and either $W_{4,1}(1,2,1)$ or
$W_{4,1}(2,1,2)$ for $p=4$. 

\begin{figure}[hb]
\includegraphics[height=1in]{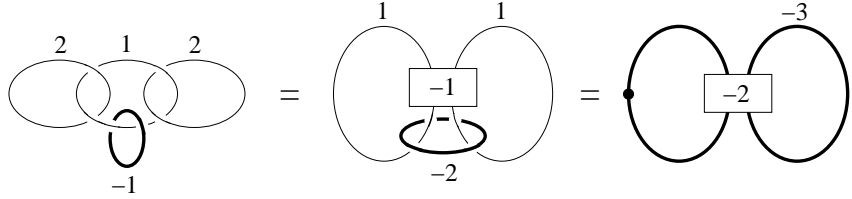} 
\caption{The manifold $W_{4,1}(2,1,2)={\mathcal C}$}
\label{f:example1} 
\end{figure}

But $W_{p,1}(1,2,\ldots,2,1)$ is diffeomorphic to the canonical
resolution $R_{p,1}=D_{-p}$, while Figure~\ref{f:example1} shows that
$W_{4,1}(2,1,2)={\mathcal C}$. As expected, we have just re--obtained
Theorem~\ref{t:mcduff-quoted}.

\begin{figure}[ht]
\includegraphics[height=1in]{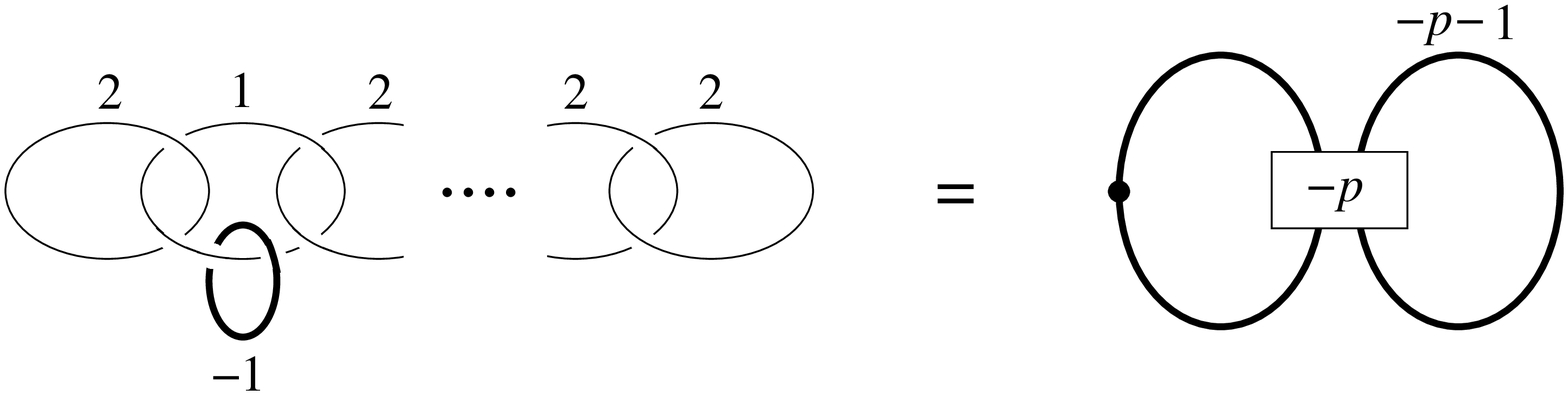} 
\caption{The manifold $W_{p^2,p-1}(p,1,2,\ldots,2)$}
\label{f:example2} 
\end{figure}

\medskip
\noindent
\textbf{Second Example.} Let us consider $L(p^2,p-1)$.
Since 
\[
\frac {p^2}{p^2-p+1} = [p,\overbrace{2,\ldots,2}^p],
\]
applying Lemma~\ref{l:0-sequence} one sees that
\[
\bZ_{p,q}=\{(1,2,\ldots,2,1), (p,1,2,\ldots,2)\}.
\]
Thus, by Theorem~\ref{t:main}, the symplectic fillings of
$(L(p^2,p-1),{\overline\xi_0})$ up to blowup and diffeomorphism are
given by the canonical resolution $R_{p^2,p-1}$ and
$W_{p^2,p-1}(p,1,2,\ldots,2)$, which is a rational homology ball as is
apparent from Figure~\ref{f:example2}. In fact,
Figure~\ref{f:example2} shows that $W_{p^2,p-1}(p,1,2,\ldots,2)$ is
precisely the rational homology ball used in the symplectic rational
blowdown construction~\cite{Sy}, and therefore carries a symplectic
form $\om$ such that $(W_{p^2,p-1}(p,1,2,\ldots,2), \om)$ is a
symplectic filling of $(L(p^2,p-1),{\overline\xi_0})$.

\begin{figure}[h]
\includegraphics[height=1in]{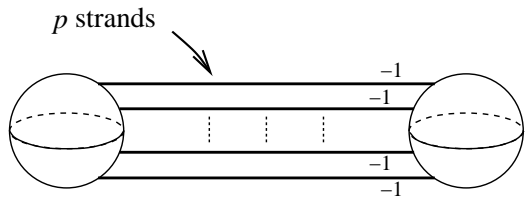} 
\caption{The manifold $W_{p,p-1}(0)$}
\label{f:example3} 
\end{figure}

\medskip
\noindent
\textbf{Third Example.} This time we look at $L(p,p-1)$. We have $k=1$
and $\bZ_{p,p-1}=\{(0)\}$. Hence, by Theorem~\ref{t:main}(b) every
filling of $L(p,p-1)$ is diffeomorphic to a blowup of $W_{p,p-1}(0)$,
which is given in Figure~\ref{f:example3}.
\smallskip
Theorem~\ref{t:main}(a) says that $W_{p,p-1}(0)$ carries a symplectic
structure $\om$ such that $(W_{p,p-1}(0),\om)$ is a symplectic filling
of $(L(p,p-1),{\overline\xi_0})$ and moreover $W_{p,p-1}(0)$ does not
contain smoothly embedded $(-1)$--spheres. In fact, $W_{p,p-1}(0)$ is
easily seen to be diffeomorphic to the canonical resolution
$R_{p,p-1}$, which gives a symplectic filling of
$(L(p,p-1),{\overline\xi_0})$ with even intersection form.

Alternatively, using~\cite{Go} one can check that $W_{p,p-1}(0)$
carries a structure of Stein surface, so in particular it carries a
symplectic structure $\om$ such that $(W_{p,p-1}(0),\om)$ is a
symplectic filling of $(L(p,p-1),\xi)$ for some contact structure
$\xi$. On the other hand, on the lens space $L(p,p-1)$ there exists
only one tight contact structure up to isotopy~\cite{Gi, Ho}. Hence,
since fillable contact structures are tight we have
$\xi={\overline\xi_0}$ and $(W_{p,p-1}(0),\om)$ is a symplectic
filling of $(L(p,p-1),{\overline\xi_0})$. As explained in the next
section, the argument we have just used can be generalized to prove
Theorem~\ref{t:main}(a).

\section{The proof of Theorem~\ref{t:main}}
\label{s:proof}

In this section we outline the proof of Theorem~\ref{t:main}. Complete
arguments will appear elsewhere~\cite{Li3}.

\medskip
\noindent\textbf{Part (a)}. The argument consists of two steps. First,
we show that each of the manifolds $W_{p,q}(\n)$ is a {\rm Stein
surface} (see e.g.~\cite{Go} for the definition).

A knot $K$ in a closed, contact three--manifold $(Y,\xi)$ is called
\emph{Legendrian} if it is everywhere tangent to the distribution
$\xi$. The contact structure induces a framing of $K$ often called the
\emph{contact framing}.

The standard tight contact structure $\ze_0$ on $S^1\x S^2$ is the
kernel of the pull--back of the one--form $zd\th + xdy - ydx$ under
the inclusion $S^1\x S^2\subset S^1\x\R^3$, where $d\th$ is the
standard angular form on $S^1$ and $x, y, z$ are coordinates on
$\R^3$. 

The manifolds $W_{p,q}(\n)$ are obtained by attaching $2$--handles to
the boundary of $S^1\x D^3$. If the $2$--handles are attached along
knots $K_i$ which are Legendrian with respect to $\ze_0$ and whose
framings equal their contact framings minus one, then for every choice
of orientations of the $K_i$'s the resulting smooth four--manifold
carries a structure of Stein surface. This construction is sometimes
referred to as {\em Legendrian surgery}. It is originally due to
Eliashberg~\cite{El1}, and has been sistematically investigated by
Gompf~\cite{Go}. The first step of the proof of
Theorem~\ref{t:main}(a) is achieved using the following
Proposition~\ref{p:legendrian}, which can be proved via an induction
argument.

\begin{prop}[\cite{Li3}]\label{p:legendrian}
Let $\widetilde\ze_0$ be the pull-back of $\ze_0$ under the
diffeomorphism ~\eqref{e:varphi} and let ${\mathbf L}$ be the
``thick'' framed link of Figure~\ref{f:W}. Then, ${\mathbf
L}$ is isotopic to a link $\LL$ Legendrian with respect to
$\widetilde\ze_0$, all of whose components have contact framing
equal to zero.
\end{prop}

The existence of a Stein structure on $W_{p,q}(\n)$ implies that
$W_{p,q}(\n)$ does not contain smoothly embedded $(-1)$--spheres.
Moreover, each Stein structure induces an exact symplectic form $\om$
on $W_{p,q}(\n)$ and a contact structure $\xi$ on $\del
W_{p,q}(\n)=L(p,q)$ such that $(W_{p,q}(\n), \om)$ is a symplectic
filling of $(L(p,q),\xi)$.

\medskip
The second step of the argument consists of showing that $W_{p,q}(\n)$
carries a Stein structure inducing a contact structure $\xi$
isomorphic to $\overline\xi_0$.

Given a $2$--plane field $\xi$ on a three--manifold $M$, there is a
$\Spin^c$ structure $\s_\xi$ on $M$ which only depends on $\xi$ up to
homotopy. The following result is proved via an induction argument
which relies on the homotopy invariant $\Ga(\xi,-)$ introduced
in~\cite{Go}.

\begin{prop}[\cite{Li3}]
The Legendrian link $\LL$ of Proposition~\ref{p:legendrian} can be
chosen and oriented so that the contact structure $\xi$ induced by 
the corresponding Stein structure satisfies $\s_\xi=\s_{\overline\xi_0}$.
\end{prop}

By the classification of tight contact structures on lens
spaces~\cite{Gi, Ho}, $\s_\xi$ determines $\xi$ up to isotopy,
hence $\xi={\overline\xi_0}$. This concludes the proof of
Theorem~\ref{t:main}(a).

\medskip
\noindent\textbf{Part (b)}. Let $(X,\om)$ be a symplectic
four--manifold. A \emph{symplectic string} in $X$ is an immersed
symplectic surface
\[
\Ga=\bigcup_{i=0}^k C_i\subset X
\]
such that:
\begin{enumerate}
\item[(a)]
$C_i$ is a connected, embedded symplectic surface, $i=0,\ldots, k$;
\item[(b)]
$C_i$ meets transversely $C_{i+1}$ at one point, $i=0,\ldots, k-1$;
\item[(c)]
$C_i\cap C_j=\emptyset$ if $|i-j|>1$.
\end{enumerate}
We say that $\Ga$ as above is \emph{of type $(m_0,\ldots,m_k)$} if,
furthermore,
\begin{enumerate}
\item[(d)]
$C_i\cdot C_i = -m_i$, $i=0,\ldots, k$.
\end{enumerate} 

\medskip
Let $(W,\om)$ be a symplectic filling of
$(L(p,q),{\overline\xi_0})$. Using a cut--and--paste argument which
combines results from~\cite{El1} and~\cite{McW}, it is possible to
prove:

\begin{thm}[\cite{Li3}]\label{t:string}
There exist a symplectic four--manifold $X_W$ and a symplectic string
\[
\Ga=\bigcup_{i=0}^k C_i\subset X_W
\]
of type $(-1,b_1-1,b_2,\ldots,b_k)$ such that $W$ is diffeomorphic to
$X_W\setminus\nu(\Ga)$, where $\nu(\Ga)\subset X_W$ is a regular
neighborhood of $\Ga$.
\end{thm}

The results of~\cite{McD1} imply that if $(M,\om)$ is a closed
symplectic four--manifold containing an embedded symplectic
$2$--sphere $C$ of self--intersection $+1$ such that $M\setminus C$ is
minimal (i.e.~not containing embedded symplectic $(-1)$--spheres),
then $(M,\om)$ is symplectomorphic to $\CP^2$ with the standard
K\"ahler structure. Morever, the symplectomorphism can be chosen so
that $C$ is sent to a complex line. Thus, since non--minimal (possibly
non--compact) symplectic four--manifolds can be reduced to minimal
ones by blowing down exceptional symplectic spheres,
applying~\cite{McD1} to the pair $(X_W, C_0)$ we conclude that, for
some $M\geq 0$, there is a symplectomorphism
\[
\psi\co X_W\to \CP^2\# M{\overline \CP}^2
\] 
such that $\psi(C_0)$ is a complex line in $\CP^2$. Clearly,
$\psi(\Ga)\subset\CP^2\# M{\overline \CP}^2$ is a symplectic string of
type $(-1,b_1-1,b_2,\ldots,b_k)$. Now we have:

\begin{thm}[\cite{Li3}]\label{t:proper-transform}
Let $\CP^2\# M\overline\CP^2$ be endowed with a blowup of the standard
K\"ahler structure on $\CP^2$. Let
\begin{equation*}
\De=\bigcup_{i=0}^k D_i \subset \CP^2\# M\overline\CP^2
\end{equation*}
be a symplectic string of type $(-1,b_1-1,b_2,\ldots,b_k)$ such that
$D_0\subset\CP^2$ is a complex line. Then,
\begin{enumerate}
\item[(a)] $\De$ is the strict transform of two distinct complex lines
in $\CP^2$;
\item[(b)] The complement of a regular neighborhood of $\De$ is
diffeomorphic to $W_{p,q}(\n)$ for some $\n\in \bZ_{p,q}$.
\end{enumerate}
\end{thm}

Note that symplectic strings have well--defined strict
transforms. Theorem~\ref{t:main}(b) follows immediately from
Theorem~\ref{t:proper-transform} together with the previous
discussion. The proof of Theorem~\ref{t:proper-transform}(a) relies on
the positivity of intersections of $J$--holomorphic
curves~\cite{McD2}, while Theorem~\ref{t:proper-transform}(b) follows
directly from Theorem~\ref{t:proper-transform}(a) via a Kirby calculus
argument.

\end{document}